\documentclass{amsart}
\usepackage{amssymb,amsmath,amsthm}
\usepackage{enumitem}
\usepackage{microtype}
\usepackage{amsfonts}
\usepackage{mathtools}

\usepackage[]{hyperref}

\usepackage{cite}

\numberwithin{equation}{section}

\newtheorem{theorem}{Theorem}[section]

\newtheorem{lemma}[theorem]{Lemma}

\theoremstyle{definition}
\newtheorem{definition}{Definition}[section]

\newcommand{\beq}{\begin{equation}}
\newcommand{\eeq}{\end{equation}}

\newcommand{\beqq}{\begin{equation*}}
\newcommand{\eeqq}{\end{equation*}}

\theoremstyle{remark}

\makeatletter
\newcommand{\Extend}[5]{\ext@arrow0099{\arrowfill@#1#2#3}{#4}{#5}}
\makeatother

\begin{document}
\hspace{1em}
\title[2D micropolar Rayleigh-B\'{e}nard convection  problem]{Global well-posedness  for  the  2D micropolar Rayleigh-B\'{e}nard convection  problem without velocity dissipation }
\author{Sheng Wang}
\address{Shanghai Center for Mathematical Sciences, Fudan University, Shanghai 200433, P.R. China;
}
\email{19110840011@fudan.edu.cn}

\begin{abstract} In this article, we study the Cauchy problem to the micropolar Rayleigh-B\'{e}nard convection  problem without velocity dissipation in two dimension. We first prove the local well-posedness of a smooth solution,  and then  establish a  blow up criterion in terms of the gradient of scalar temperature field. At last, we obtain the global well-posedness to the system.
\\[0.6em]
\\[0.6em]

\end{abstract}
\maketitle
\section{introduction}\label{section-intro}

 \setcounter{section}{1}\setcounter{equation}{0}
In fluid dynamics, the micropolar fluid model can be regarded as an extension of the classical fluid dynamics model represented by the Navier-Stokes equation, because the classical fluid dynamics model can't describe the microstructure. In 1965, C.A. Eringen first introduced equations to simulate micropolar fluids (see\cite{ACE}).
Interestingly,  many experiments have shown that this model is more accurate than classical fluid dynamics models in describing some fluid movements, such as blood (see\cite{PRU}). When we investigate the behavior of a fluid layer filling the area between two rigid surfaces,  the fluid is heated from below. We usually describe this heat convection phenomenon in the framework of the Boussinesq equation approximation,  which is called micropolar Rayleigh-B\'{e}nard convection  problem(see\cite{FMT,KLL,TA}) and  described by the following equations  in $\mathbb{R}^{3}$



\begin{equation*}
\left\{
\begin{aligned}
&u_t + u\cdot \nabla u - (\nu+\kappa)\Delta u + \nabla P = 2\kappa \nabla\times \omega + e_{3}\cdot\theta,\\
&\omega_t + u\cdot \nabla \omega - \gamma\Delta \omega-(\alpha+\beta)\nabla div \omega + 4\kappa\omega = 2\kappa \nabla\times\, u, \\
& \theta_t + u\cdot \nabla \theta - \mu\Delta \theta  = u\cdot e_{3}, \\
&\operatorname{div} u=0,\\
&(u,\,\omega,\,\theta)(x,t)_{t=0}=(u_{0},\,\omega_{0},\,\theta_{0})(x) \quad x\in \mathbb{R}^3.
\end{aligned}
\right.
\end{equation*}
where $u$ is the fluid velocity, $\omega$ is the field of microrotation representing the angular velocity of the rotation of the particles of the fluid, $P$ is the scalar pressure of the flow, $\theta$ is the scalar temperature, $\nu$ is the Newtonian kinematic viscosity, $\kappa$ is the micro-rotation viscosity, $\alpha$ and $\beta$
are the angular viscosities, $\mu$ is the thermal diffusivity, $e_3= (0,0,1)$ is the vertical unit vector, the forcing term $e_3\cdot\theta$ in the momentum equation describes the action of the buoyancy force on fluid motion, and $u\cdot e_3$  models the Rayleigh-B\'{e}nard convection in a heated fluid.
Many recent efforts are focused on the 2D micropolar fluid with partial
dissipation \cite{DZ,DLW,XUE}.

In present article, we shall consider  the 2D micropolar Rayleigh-B\'{e}nard convection  problem. For convenience, we assume that the velocity component in the $x_{3}$-direction is
zero and the axes of rotation of particles are parallel to the $x_{3}$-axis. That is,
$$u=u(x_{1},x_{2},0),\quad \omega=\omega(0,0,\omega_{3}(x_{1},x_{2})),\quad P=P(x_{1},x_{2}), \quad \theta=\theta(x_{1},x_{2}),$$
which implies that $$\operatorname{div} \omega=0,\quad \nabla\times u= \partial_{1}u_{2}-\partial_{2}u_{1}, \quad \nabla\times \omega=(\partial_{2}\omega, -\partial_{1}\omega), \quad \nabla\times\nabla\times \omega=-\Delta\omega.$$
Thus, the 2D micropolar Rayleigh-B\'{e}nard convection  problem can be written as
\begin{equation}\label{equ:1.1-1}
\left\{
\begin{aligned}
&u_t + u\cdot \nabla u - (\nu+\kappa)\Delta u + \nabla P = 2\kappa \nabla\times \omega + e_{2}\cdot\theta,\\
&\omega_t + u\cdot \nabla \omega - \gamma\Delta \omega + 4\kappa\omega = 2\kappa \nabla\times\, u, \\
& \theta_t + u\cdot \nabla \theta - \mu\Delta \theta  = u\cdot e_{2}, \\
&\operatorname{div} u=0,\\
&(u,\,\omega,\,\theta)(x,t)_{t=0}=(u_{0},\,\omega_{0},\,\theta_{0})(x).
\end{aligned}
\right.
\end{equation}

In 2020, Xu et al. \cite{XU1} first obtained global regularity for the  system \eqref{equ:1.1-1}  with zero diffusivity (i.e., $\mu=0$).
In 2021,  Xu et al. \cite{XU2} proved the unique local solvability  of smooth solution  when the system \eqref{equ:1.1-1} has only  velocity  dissipation,   and then established a criterion for the breakdown of smooth solutions imposed only the maximum norm of the gradient of scalar temperature field. Furthermore,  they finally showed the global regularity of the system  with zero angular viscosity (i.e., $\gamma=0$).

However,  the related results  is still unknown  for the case $\nu+\kappa=0$.
The purpose of this work are to study  the well-posedness of the  following system \eqref{equ:1.1-1} without velocity dissipation.

\begin{equation}\label{equ:1.1}
\left\{
\begin{aligned}
&u_t + u\cdot \nabla u + \nabla P = 2\kappa \nabla\times \omega + e_{2}\cdot\theta,\\
&\omega_t + u\cdot \nabla \omega -\gamma\Delta\omega + 4\kappa\omega = 2\kappa \nabla\times\, u, \\
& \theta_t + u\cdot \nabla \theta -\mu\Delta\theta = u\cdot e_{2}, \\
&\operatorname{div} u=0,\\
&(u,\,\omega,\,\theta)(x,t)_{t=0}=(u_{0},\,\omega_{0},\,\theta_{0})(x).
\end{aligned}
\right.
\end{equation}

We first prove the local well-posedness of the smooth solution for  the system \eqref{equ:1.1}, and then  establish a blow up  criterion in terms of the gradient of scalar temperature field and finally obtain  global well-posedness of the system.

The main results in this paper are stated as follows.
\begin{theorem}\label{th1}(Local well-posedness) Let $\gamma>0, \mu>0$,  and $\operatorname{div} u_0=0$. Suppose  $s>2$,  and $(u_0, \omega_0, \theta_0,
	)\in H^s(\mathbb{R}^2)$. Then there exist  $T>0$ and a unique
	solution $(u,\omega,\theta)$ of the system \eqref{equ:1.1} such that
	\begin{equation*}
	u\in C\big([0, T]; H^s(\mathbb{R}^2)\big),\quad
	(\omega,\theta)\in C\big([0, T]; H^s(\mathbb{R}^2)\big)\cap L^{2}\big(0, T;
	H^{s+1}(\mathbb{R}^2)\big).
	\end{equation*}
\end{theorem}
\begin{theorem}\label{th2}(Blow-up criterion ) Let $\gamma>0, \mu>0$, and $\operatorname{div} u_0=0$. Suppose $s>2$,  $(u_0, \omega_0,\theta_0,
	)\in H^s(\mathbb{R}^2)$  and solution $(u,\omega,\theta)$ of the system \eqref{equ:1.1} in $ H^s(\mathbb{R}^2)$ satisfies
	\begin{equation}\label{blow-up}
	\|\nabla \theta(.,\tau)\|_{L^{\infty}}<+\infty,
	\end{equation}
	then $(u,\omega,\theta)\in H^s(\mathbb{R}^2) $ for all time $t\in [0; T]$.
\end{theorem}
\begin{theorem}\label{th3}(Global regularity) Let $\gamma>0, \mu>0$ and $\operatorname{div} u_0=0$. Suppose   $s>2$,  and $(u_0, \omega_0,\theta_0,
	)\in H^s(\mathbb{R}^2)$. Then the system \eqref{equ:1.1}  has a unique global
	solution $(u,\omega,\theta)$  such that for any $T>0$,
	\begin{equation*}u\in C\big([0, \infty); H^s(\mathbb{R}^2)\big),\quad
	(\theta,\omega) \in C\big([0, \infty); H^s(\mathbb{R}^2)\big)\cap L^{2}\big(0, T;
	H^{s+1}(\mathbb{R}^2)\big).\end{equation*}
\end{theorem}
Here, let us explain some of the main difficulties and techniques involved in the process. First, the estimate of $\|\Omega\|_{L^{p}}:=\|\nabla \times u\|_{L^{p}}$$(2\leq p \leq\infty)$ plays a very key role in our proof. In order to get it, the main  difficulty is due to the dynamic micro-rotational term $\nabla\times \omega$ in the velocity equation. Indeed, taking a $\nabla\times$ on velocity   equation of \eqref{equ:1.1},  we find that
$$\Omega_t + u\cdot\nabla \Omega + 2\kappa \Delta  \omega =\partial_1\theta.$$
However, due to the non-dissipation property of the variable $u$, we have no obvious way to handle the term  $\kappa \Delta  \omega$ as a perturbation and force us to deal with $2\kappa \Delta  \omega$  as a whole.
Here, to obtain the estimate of  $\|\Omega\|_{L^{\infty}}$ so that we can get the estimate of $\|\nabla u\|_{L^{\infty}}$ and show the global regularity of the 2D micropolar fluid flows without velocity disspation, we used the idea in \cite{DLW} and introduced a nice quantity $Z=\Omega+\frac{2\kappa}{\gamma}\omega$ to  overcome the difficulty.
Fortunately,  the key observation for the more complicated system  \eqref{equ:1.1} is still effective. Resort to  the quantity $Z$, we  can exploit the estimates of $\|\Omega\|_{L^{\infty}}$ and $\|\Omega\|_{L^{p}}$ and then establish a criterion for the breakdown of smooth solutions imposed only the maximum norm of the gradient of scalar temperature field by the classical  logarithm inequality. In particular, this demonstrates that the extension of  smooth solutions system  \eqref{equ:1.1} is dominated by the temperature field $\theta$ of the fluid. In the proof of the global regularity of the system  without velocity dissipation,  using the heat kernel estimation,
we can obtain $\|\Omega\|_{L^p}\leq C(u_0,\omega_0,\theta_0;T)$$(2\leq p<\infty)$. Furthermore, combining  with Gagliardo-Nirenberg interpolation and the Calderon-Zygmund inequalities, one can obtain the  key estimate of  $\|u\|_{L^\infty}$. Based on it, making  full use of the dissipation property of the variable $\theta$, we finally get the desired estimate $\|\nabla\theta \|_{L^\infty}$.

The rest of this article follows. First, we will  recall some facts about the Littlewood-Paley decomposition, the frequency characterization of Sobolev space and Besov space, and the estimation of some nonlinear terms.
Then, in Section 3,  we present the  proof of Theorem \ref{th1}. Section 4 is devoted to the proof of Theorem \ref{th2}. In the last section, we shall prove  the global regularity of the system without velocity dissipation.

\par
\section{ Preliminaries and Lemmas  }

In the preparatory  section, we introduce some
common notations,  some basic points about Littlewood-Paley theory
and compile some auxiliary lemmas.

Let $\mathcal {S}(\mathbb{R}^{2})$ be the Schwartz class of rapidly
decreasing function. Given $f\in \mathcal {S}(\mathbb{R}^{2})$, its
Fourier transform $\mathcal {F}f= \widehat{f}$ is defined by
$$\widehat{f}(\xi)=\int_{\mathbb{R}^{2}}e^{-ix\cdot\xi}f(x)dx.$$
Let $(\chi, \varphi)$ be a couple of smooth functions valued in $[0,
1]$ such that $\chi$ is supported in the ball $\{\xi\in
\mathbb{R}^{2}:  \ |\xi|\leq\frac{4}{3}\}$, $\varphi$ is supported
in the shell $\{\xi\in \mathbb{R}^{2}: \ \frac{3}{4}
\leq|\xi|\leq\frac{8}{3}\}$,  $\varphi(\xi)=\chi(\xi/2)-\chi(\xi)$
and
$$\chi(\xi)+\sum_{j\geq0}\varphi(2^{-j}\xi)=1, \ \forall \, \xi \in\mathbb{ R}^{2},$$
$$\sum_{j\in Z}\varphi(2^{-j}\xi)=1, \ \forall \, \xi \in \mathbb{R}^{2}\setminus\{0\}.$$

 First, we introduce the dyadic blocks $\Delta_j$  by setting
\begin{equation*}
\Delta_j f=0\ \hbox{ if }\ j\leq-2,\quad \Delta_{-1}f={\mathcal
	F}^{-1}(\chi{\mathcal F}f) \quad\hbox{and}\quad \Delta_j f={\mathcal
	F}^{-1}(\varphi(2^{-j}\cdot){\mathcal F} f)\ \text{ if } \  j\geq 0.
\end{equation*}
One may prove  that for all tempered distribution $u$ the following
Littlewood-Paley decomposition holds true:
\begin{equation*}f=\sum\limits_{j\geq -1}\Delta_j f, \quad
{S}_{j}f=\sum\limits_{-1\leq q\leq j-1}\Delta_q f.
\end{equation*}
One easily verifies that with our choice of $\varphi$,
\begin{equation*}
\Delta_{j}\Delta_{q}f\equiv0,\quad \hbox{if} \quad |j-q|\geq2\quad
and \quad \Delta_{j}(S_{q-1}f\Delta_{q}f)\equiv0,\quad \hbox{if}
\quad |j-q|\geq5.
\end{equation*}
Let us recall the definition of the Besov space.
\begin{definition} Let $s\in \mathbb{R}$, $1\le p,
	r\le+\infty$. The nonhomogeneous Besov space $B^{s}_{p,r}$ is
	defined by
	$$B^{s}_{p,r}=\Big\{f\in {\mathcal {S}'}(\mathbb{R}^2):\,\|f\|_{B^{s}_{p,r}}<+\infty\Big\},$$
	where $$ \|f\|_{B^{s}_{p,r}}:= \Big\|2^{js}
	\|\Delta_jf\|_{p}\Big\|_{\ell^r}.$$
\end{definition}

Notice that the usual Sobolev space $H^s$ coincides with $B_{2,
	2}^s$ for every $s\in\mathbb{R},$ and the norm of $H^s$ can be characterized by
\begin{equation*}\|f\|_{H^s}\approx\|{S}_{0}f\|_{L^{2}} +\Big(\sum_{q\geq0}2^{2qs}\|\Delta_q f\|^{2}_{L^{2}} \Big)^{\frac{1}{2}}.
\end{equation*}
We shall also use  the following some important lemmas in the
proof of our main results.

\begin{lemma}\label{l:commutator} We present some commutator estimates as follows
	
	\noindent{(i)}(Lemma 2.2 of \cite{DZ})\ Let $s>0$, $f\in H^{s}(\mathbb{R}^{2})\cap
	L^{\infty}(\mathbb{R}^{2})$, and $g\in H^{s+1}(\mathbb{R}^{2})\cap
	W^{1,\infty}(\mathbb{R}^{2})$ with $\operatorname{div} g=0$. Then the following calculus
	inequality holds $$\|[\Delta_j, g]\cdot\nabla f\|_{L^{2}}\leq
	Cc_j2^{-js}\Big(\|\nabla
	g\|_{L^{\infty}}\|f\|_{H^{s}}+\|\nabla g\|_{H^{s}}\|f\|_{L^{\infty}}\Big),$$
	where $c_j$ is a sequence satisfying $\|{c_j}\|_{\ell^{2}}\leq1$.
	
	\noindent{(ii)}((Lemma 2.100 of \cite{BCD})\ Let $s>0$, $f\in H^{s}(\mathbb{R}^{2})\cap
	W^{1,\infty}(\mathbb{R}^{2})$, and $g\in H^{s}(\mathbb{R}^{2})\cap
	W^{1,\infty}(\mathbb{R}^{2})$. Then the following calculus
	inequality holds $$\|[\Delta_j, g]\cdot\nabla f\|_{L^{2}}\leq
	Cc_j2^{-js}\Big(\|\nabla
	g\|_{L^{\infty}}\|f\|_{H^{s}}+\| g\|_{H^{s}}\|\nabla f\|_{L^{\infty}}\Big),$$
	where $c_j$ is a sequence satisfying $\|{c_j}\|_{\ell^{2}}\leq1$.
\end{lemma}
\begin{lemma}(\cite{Che})\label{bernstein} Let $j, k\in\mathbb{N},$ $1\leq a\leq
	b$ and $f\in L^a(\mathbb{R}^2).$
	There exists a constant $C$ such that the  following  inequalities hold
	\begin{eqnarray*}
		\sup_{|\alpha|=k}\|\partial ^{\alpha}S_{j}f\|_{L^b}&\leq& C^{k}\,2^{j(k+2(\frac{1}{a}-\frac{1}{b}))}\|S_{j}f\|_{L^a},\\
		\ C^{-k}2^
		{jk}\|{\Delta}_{j}f\|_{L^a}&\leq&\sup_{|\alpha|=k}\|\partial
		^{\alpha}{\Delta}_{j}f\|_{L^a}\leq
		C^{k}2^{jk}\|{\Delta}_{j}f\|_{L^a}.
	\end{eqnarray*}
\end{lemma}
\begin{lemma}\label{morse}(\cite{KG})\ Let $s>0$, $f, g\in H^{s}(\mathbb{R}^{2})\cap
	L^{\infty}(\mathbb{R}^{2})$. Then the  following  inequality holds
	$$\|fg\|_{H^{s}}\leq
	C\Big(\|f\|_{L^{\infty}}\|g\|_{H^{s}}+\|f\|_{H^{s}}\|g\|_{L^{\infty}}\Big).$$
\end{lemma}
To prove this proposition, we recall the maximal regularity property for the heat kernel (see, \cite{HM}).
\begin{lemma}
	The operator T defined by
	\\$$Tf(t)=\int^t_{0}\Delta e^{(t-s)\Delta} f(s)ds$$\\
	maps $L^p(0,T;L^q(\mathbb{R}^d))$ to $L^p(0,T;L^q(\mathbb{R}^d))$ for any T$\in(0,\infty)$ and $p,q\in(1,\infty)$.
\end{lemma}

\section{The  proof of Theorem \ref{th1}}
\ \ \ \ \  This section is devoted to the proof of the local well-posedness  for a smooth solution to the system \eqref{equ:1.1}.  For a clear presentation, we
split it into the following several steps.

\noindent{\bf Step 1.}\  Construction of an approximate solution
sequence.

Here we use the classical Friedrich's method: For $n\geq 1$, let
$J_n$ be the spectral cut-off defined by
\begin{equation*}
\widehat{J_{n}f}(\xi)=1_{[0,n]}(|\xi|)\widehat{f}(\xi), \quad \xi
\in\mathbb{ R}^2.
\end{equation*}
We consider the following system in the space $L^{2}_{n}:=\{f\in
L^2(\mathbb{R}^{2})|\text{ supp} f\subset B(0,n)\}$:
\begin{align}\label{eq 3.1}
\left\{
\begin{aligned}
&\partial_tu_n+\mathcal{P}J_{n}(J_nu_n\cdot \nabla J_nu_n)=2\kappa \mathcal{P}\nabla\times \omega_n+\mathcal{P}J_n(\theta_n\cdot e_2),\\
&\partial_t\omega_n +J_n(J_nu_n \cdot\nabla J_n\omega_n)-\gamma\Delta \omega_n +4\kappa\omega_{n}=2\kappa \mathcal{P}\nabla\times u_n,\\
&\partial_t\theta_n +J_n(J_nu_n \cdot\nabla J_n\theta_n)-\mu\Delta\theta_n=\mathcal{P}J_n(u_n \cdot e_2),\\
&(u_n,\omega_n,\theta_n)|_{t=0}=J_{n}(u_0,\omega_0,\theta_0).
\end{aligned}
\right.
\end{align}
Here $\mathcal{P}$ denotes the Helmholtz projection operator. The
Cauchy-Lipschitz theorem entails that the system \eqref{eq 3.1}
exists a unique maximal solution $(u_n,\omega_n,\theta_n)$ in
$\mathcal{C}^{1}([0,T^{*}_{n});L^{2}_{n})$. On the other hand, we
observe that $J^{2}_n=J_n, \mathcal{P}^{2}=\mathcal{P}$ and
$J_n\mathcal{P}=\mathcal{P}J_n$. It follows that
$(\mathcal{P}u_n,\omega_n,\theta_n)$ and
$(J_n\mathcal{P}u_n,J_n\omega_n,J_n\theta_n)$ are also solutions to
the system \eqref{eq 3.1}. The uniqueness gives that $ J_nu_n=u_n, J_n\omega_n=\omega_n $ and $ J_n\theta_n=\theta_n$. Therefore,
\begin{align}\label{eq 3.2}
\left\{
\begin{aligned}
&\partial_tu_n+\mathcal{P}J_{n}(u_n\cdot \nabla u_n)=2\kappa \mathcal{P}\nabla\times \omega_n+\mathcal{P}J_n(\theta_n\cdot e_2),\\
&\partial_t\omega_n +J_n(u_n \cdot\nabla \omega_n)-\gamma \Delta \omega_n+4\kappa\omega_{n}=2\kappa \mathcal{P}\nabla\times u_n,\\
&\partial_t\theta_n +J_n(u_n \cdot\nabla \theta_n)-\mu\Delta\theta_n=\mathcal{P}J_n(u_n\cdot e_2),\\
&(u_n,\omega_n,\theta_n)|_{t=0}=J_{n}(u_0,\omega_0,\theta_0).
\end{aligned}
\right.
\end{align}
As the operators $J_n$ and $\mathcal{P}J_n$ are the orthogonal
projectors for the $L^2$-inner product, the above formal
calculations remain unchanged.

\noindent{\bf Step 2.} Energy estimates.

First, we  easily obtain the following basic  $L^{2}$-estimates of  the system \eqref{eq 3.2}
\begin{equation}\label{eq 3.3}
\begin{array}{ll}
\left\Vert u_n\right\Vert_{L^{2}}^{2}+\left\Vert \omega_n\right\Vert_{L^{2}}^{2}+\left\Vert \theta_n\right\Vert_{L^{2}}^{2}+2\gamma\displaystyle\int_{0}^{t}\|\nabla \omega_n(\tau)\|^{2}_{L^{2}}d\tau+2\mu\displaystyle\int_{0}^{t}\|\nabla \theta_n(\tau)\|^{2}_{L^{2}}d\tau
\leq C(u_0,\omega_0,\theta_0;T).
\end{array}
\end{equation}

Next we present the uniform estimates for the approximate solutions $ (u_n,\omega_n,\theta_n)$ in $H^{s}$. Applying the operator $\Delta_j(j\geq0)$ to the
system \eqref{eq 3.2}, and taking the $L^2$-scalar product of the
first equation of \eqref{eq 3.2} with $\Delta_ju_{n}$,
the second equation with $\Delta_j\omega_{n}$, and the third equation with  $\Delta_j\theta_{n}$ respectively, by the
condition $\text{div}u_{n}=0$, we obtain
\begin{equation*}
\label{3.16}\begin{split}&\frac{1}{2}\frac{d}{dt}\big(\|\Delta_ju_{n}\|_{L^{2}}^{2}
+\|\Delta_j\omega_{n}\|_{L^{2}}^{2}+\|\Delta_j\theta_{n}\|_{L^{2}}^{2}\big)+ \gamma\|\nabla\Delta_j\omega_{n}\|_{L^{2}}^{2}+\mu\|\nabla\Delta_j\theta_{n}\|_{L^{2}}^{2}\\
&=-\big([\Delta_j, u_n] \nabla u_n, \Delta_ju_n
\big)_{L^{2}}-\big([\Delta_j, u_n] \nabla \omega_n, \Delta_j\omega_n
\big)_{L^{2}}-\big([\Delta_j, u_n] \nabla \theta_n, \Delta_j\theta_n
\big)_{L^{2}}\\
&\quad+\kappa\big(\Delta_j\nabla\times \omega_n, \Delta_ju_n
\big)_{L^{2}}+\kappa\big(\Delta_j\nabla\times u_n, \Delta_j\omega_n
\big)_{L^{2}}\\&\quad+\big(\Delta_j\theta_n, \Delta_ju_n
\big)_{L^{2}} +\big(\Delta_ju_n, \Delta_j\theta_n
\big)_{L^{2}}.
\end{split}
\end{equation*}
Employing Lemma \ref{l:commutator}  and  the Young
inequality, multiplying both sides by $2^{2js}$, summing over $j\geq 0$, and then combining with  \eqref{eq 3.3}, we thus get
\begin{equation}
\label{eq 3.4}\begin{split}&\frac{1}{2}\frac{d}{dt}\big(\|u_n\|_{H^{s}}^{2}
+\|\omega_n\|_{H^{s}}^{2}+\|\theta_n\|_{H^{s}}^{2}\big)+ \gamma\|\nabla \omega_{n}\|_{H^{s}}^{2}+\mu\|\nabla \theta_{n}\|_{H^{s}}^{2}
\\&\leq C \big(1+ \|\nabla u_n\|_{L^{\infty}}+\|\nabla\omega_n\|_{L^{\infty}}+\|\nabla
\theta_n\|_{L^{\infty}}\big)\big(\|u_n\|_{H^{s}}^{2}
+\|\omega_n\|_{H^{s}}^{2}+\|\theta_n\|_{H^{s}}^{2}\big),
\end{split}
\end{equation}
 from which
and  the Gronwall inequality, it follows that
\begin{equation}\label{eq 3.5}
\begin{array}{lll}E_n(t)\leq \big(\|u_0\|_{H^{s}}^{2}
+\|\omega_0\|_{H^{s}}^{2}+\|\theta_0\|_{H^{s}}^{2}\big)\exp\Big(\displaystyle\int_{0}^{t}G_{n}(\tau)d\tau\Big),\end{array}\end{equation}
where \begin{align*} &E_n(t)=\big(\|u_n\|_{H^{s}}^{2}
+\|\omega_n\|_{H^{s}}^{2}+\|\theta_n\|_{H^{s}}^{2}\big)+2\gamma\displaystyle\int_{0}^{t}\|\nabla \omega_{n}(\tau)\|_{H^{s}}^{2}d\tau+2\mu\displaystyle\int_{0}^{t}\|\nabla\theta _{n}(\tau)\|_{H^{s}}^{2}d\tau,\\
&G_{n}(t)=1+ \|\nabla u_n(t)\|_{L^{\infty}}+\|
\nabla\omega_n(t)\|_{L^{\infty}}+\|\nabla
\theta_n(t)\|_{L^{\infty}}.\end{align*}
\noindent{\bf Step 3.} Uniform estimates and existence of the
solution.
We denote $T_{n}^{\ast}$  by the maximal existence time of the
solution $ (u_n, \omega_n, \theta_n)$ and define
$$\tilde{T}_{n}^{\ast}=\sup\Big\{t\in [0, T_{n}^{\ast}): E_n(\tau)\leq2\big(\|u_0\|_{H^{s}}^{2}
+\|\omega_0\|_{H^{s}}^{2}+\|\theta_0\|_{H^{s}}^{2}\big)\quad\hbox{for}\quad
\tau\in[0,t]\Big \}.$$ From \eqref{eq 3.5} and  the Sobolev embedding ($s>2$), we
find that
$$E_n(t)\leq \big(\|u_0\|_{H^{s}}^{2}
+\|\omega_0\|_{H^{s}}^{2}+\|\theta_0\|_{H^{s}}^{2}\big)\exp\Big(C\big(1+\|u_0\|_{H^{s}}^{2}
+\|\omega_0\|_{H^{s}}^{2}+\|\theta_0\|_{H^{s}}^{2}\big)t\Big), t\in
[0,\tilde{T}_{n}^{\ast}).$$  Taking $T$  small enough such that
$$\exp\big(C(1+\|u_0\|_{H^{s}}^{2}
+\|\omega_0\|_{H^{s}}^{2}+\|\theta_0\|_{H^{s}}^{2})T\big)\leq\frac{3}{2}.$$
Now we can conclude that $\tilde{T}_{n}^{\ast}\geq T$. Otherwise, we
have
$$E_n(t)\leq\frac{3}{2} \big(\|u_0\|_{H^{s}}^{2}
+\|\omega_0\|_{H^{s}}^{2}+\|\theta_0\|_{H^{s}}^{2}\big)\quad\hbox{for}\quad
t\in [0,\tilde{T}_{n}^{\ast}],$$ which contradicts with the
definition of $\tilde{T}_{n}^{\ast}$. Thus, the approximate solution
$ (u_n, \omega_n, \theta_n)$ exists on $[0,T ]$ and satisfies the
following uniform estimates
\begin{equation}\label{eq 3.6}
E_n(t)\leq2 \big(\|u_0\|_{H^{s}}^{2}
+\|\omega_0\|_{H^{s}}^{2}+\|\theta_0\|_{H^{s}}^{2}\big)\quad\hbox{for}\quad
t\in [0,T].\end{equation} On the other hand, it is easy to verify by
using equations \eqref{eq 3.2},  uniform estimates \eqref{eq 3.6}
and Lemma \ref{morse} that $(\partial_{t} u_n,\partial_{t}
\omega_n,\partial_{t} \theta_n )$ is uniformly bounded in $L^{2}(0,T;
H^{s-1}(\mathbb{R}^{2}))$. Thus, the Aubin-Lions compactness theorem
(see e.g. \cite{Aubin}) ensures that there exists a subsequence $(
u_{n_{k}}, \omega_{n_{k}}, \theta_{n_{k}} )$ of $( u_n,\omega_n,
\theta_n )$ converging to a family function $( u,\omega,\theta)$ such
that 
$$u\in L^{\infty}(0, T;\\ H^s(\mathbb{R}^2)), (\omega,\theta)\in
L^{\infty}(0, T; H^s(\mathbb{R}^2))\cap L^{2}(0, T;
H^{s+1}(\mathbb{R}^2)).$$  Then passing to limit in \eqref{eq 3.2}, it
is easy to see that $(u, \omega,\theta)$ satisfies  the system
\eqref{equ:1.1}.

\noindent{\bf Step 4.} Continuity in time of the solution.

Revisiting the proof of Step 2, one can in fact obtain better
uniform estimates for $(u_n, \omega_n, \theta_n)$ (thus for $(u,
\omega, \theta)$):
\begin{equation*}\|u\|_{ \tilde{L}^{\infty}(0, T;  H^s(\mathbb{R}^2))}+\|\omega\|_{ \tilde{L}^{\infty}(0, T;  H^s(\mathbb{R}^2))}+\|\theta\|_{ \tilde{L}^{\infty}
	(0, T;  H^s(\mathbb{R}^2))}\leq C,\end{equation*} where $\|f\|_{
	\tilde{L}^{\infty}(0, T; H^s(\mathbb{R}^2))}=
\sum_{j\geq-1}2^{2j
	s}\|\Delta_{j}f\|^{2}_{L^{\infty}(0, T;L^{2}) }.$
Then we can conclude $(u, \omega,\theta)\in C([0,T ]; H^{s}
(\mathbb{R}^{2}))$. Indeed, for any $\varepsilon >0$, we can take
$N$ big enough such that $$\sum_{j>
	N}2^{2js}\|\Delta_{j}u\|^{2}_{L^{\infty}(0, T;L^{2})}\leq
\frac{\varepsilon}{4}.$$ For any $t\in (0,T )$ and $\sigma$ such
that $t +\sigma \in [0,T ]$, we have
\begin{equation}\label{eq 3.7}
\begin{array}{lll}\|u(t+\sigma)-u(t)\|_{H^{s}}^{2}&\leq \sum_{-1\leq j\leq
	N}2^{2js}\|\Delta_{j}\big(u(t+\sigma)-u(t)\big)\|^{2}_{L^{2}}+\displaystyle\frac{\varepsilon}{2}
\cr\noalign{\vskip 3mm}&\leq\sum_{-1\leq j\leq
	N}2^{2js}|\sigma|\|\partial_{t}u\|^{2}_{L^{2}(0,T;
	L^{2})}+\displaystyle\frac{\varepsilon}{2} \cr\noalign{\vskip
	3mm}&\leq2N2^{2Ns}|\sigma|\|\partial_{t}u\|^{2}_{L^{2}(0,T;
	L^{2})}+\displaystyle\frac{\varepsilon}{2}\leq\varepsilon
,\end{array}\end{equation} for $|\sigma|$ small enough. Hence,
$u(t)$ is continuous in $H^{s} (\mathbb{R}^{2})$ at the time $t$.
Similarly, we can obtain $(\omega,\theta)\in C([0, T];
H^s(\mathbb{R}^{2}))$.

\noindent{\bf Step 5.} Uniqueness of the solution.

Let $(u_1, \omega_1, \theta_1)$ and $(u_2, \omega_2, \theta_2)$ be two
solutions of \eqref{equ:1.1} with the same initial data. We denote
$\delta u = u_1-u_2,$ $\delta\omega
= \omega_1 -\omega_2, \delta\theta = \theta_1 -\theta_2$  and $\delta P = 0$. Then $(\delta u,
\delta\omega, \delta\theta )$ satisfies
\begin{equation}\label{eq3.8}
\left\{\begin{array}{lll}
\partial_t\delta u+u_{2}\cdot\nabla \delta u=-\delta u\cdot\nabla u_1+2\kappa\nabla\times\delta \omega+\delta\theta e_2,\cr\noalign{\vskip 2mm}
\partial_t\delta\omega +u_{2}\cdot\nabla \delta \omega-\gamma \Delta \delta\omega+4\kappa\delta\omega=-\delta u\cdot\nabla \omega_1+2\kappa\nabla\times\delta u,\cr\noalign{\vskip 2mm}
\partial_t\delta\theta +u_{2}\cdot\nabla \delta \theta-\mu\Delta\delta\theta=-\delta u\cdot\nabla \theta_1+\delta u e_2,\cr\noalign{\vskip 2mm}
(\delta u, \delta\theta, \delta\phi )|_{t=0}=(0,0,0).
\end{array}\right.\end{equation}

Taking $L^{2}(\mathbb{R}^{2})$ energy estimates, it is easy to show
that
$$\frac{d}{dt}(\|\delta u\|^{2}_{L^{2}}+\|\delta \omega\|^{2}_{L^{2}}+
\|\delta \theta\|^{2}_{L^{2}})\leq C (\|\delta
u\|^{2}_{L^{2}}+\|\delta \omega\|^{2}_{L^{2}}+\|\delta
\theta\|^{2}_{L^{2}}),$$ which along with the  Gronwall inequality implies
$(\delta u, \delta\omega, \delta\theta )=(0,0,0)$.
The proof of Theorem \ref{th1} is completed.

\section{The  proof of Theorem \ref{th2}}
This section is devoted to the proof of
blow-up criterion for smooth solutions to the system \eqref{equ:1.1}.
First, for the basic $L^2$-estimate  and the $\dot{H^1}$-estimate  of the system (1.1), we readily get

\begin{equation}\label{eq4.1}
\begin{array}{ll}
\left\Vert u\right\Vert_{L^{2}}^{2}+\left\Vert \omega\right\Vert_{L^{2}}^{2}+\left\Vert \theta\right\Vert_{L^{2}}^{2}+\displaystyle\int_{0}^{t}\big(\left\Vert \nabla \omega \right\Vert_{L^{2}}^{2}
+\left\Vert \nabla \theta\right\Vert_{L^{2}}^{2}\big)d\tau\leq C(u_0,\omega_0,\theta_0;T),
\end{array}
\end{equation}
and

\begin{equation}\label{eq4.2}
\begin{array}{ll}
\left\Vert \nabla u\right\Vert_{L^{2}}^{2}+\left\Vert \nabla \omega\right\Vert_{L^{2}}^{2}+\left\Vert \nabla \theta\right\Vert_{L^{2}}^{2}+\displaystyle\int_{0}^{t}\big(\left\Vert \nabla^{2} \omega \right\Vert_{L^{2}}^{2}
+\left\Vert \nabla^{2} \theta\right\Vert_{L^{2}}^{2}\big)d\tau\leq C(u_0,\omega_0,\theta_0;T).
\end{array}
\end{equation}
Denote the vorticity of velocity $u$ by $\Omega=\nabla\times u$. Applying the operator $\nabla\times$  to the first equation of \eqref{equ:1.1} yield that
$$\Omega_t + u\cdot\nabla \Omega +2\kappa \Delta  \omega =\partial_1\theta.$$ Due to the lack of volecity viscosity, we have no obvious way to handle
$ \Delta  \omega$. As in \cite{DLW}, set
\begin{equation}\label{eq3.2}Z=\Omega+\frac{2\kappa}{\gamma}\omega.\end{equation}
Thus, we can take advantage of its evolution in time governed by
\begin{equation}\label{eq3.3}Z_t + u\cdot\nabla Z =\frac{4\kappa^{2}}{\gamma}Z-\frac{16\kappa^{3}}{\gamma}\omega+\partial_1\theta.\end{equation}
This leads to the following $L^{p}$-estimate of $Z$ after multiplying \eqref{eq3.3} by $|Z|^{p-2}Z$ with $p\geq2$, integrating over $\mathbb{R}^{2}$, using the incompressibility of $u$ and  the H\"{o}lder inequality, we have

\begin{equation*}
\label{3.25}\begin{split}
&\frac{1}{p}\frac{d}{dt}\|Z\|^p_{L^p} \leq \frac{16\kappa^{3}}{\gamma}\|\omega\|_{L^p}\|Z\|^{p-1}_{L^p}+\frac{4\kappa^{2}}{\gamma}\|Z\|^{p}_{L^p}+\|\nabla\theta\|_{L^p}\|Z\|^{p-1}_{L^p}.
\end{split}
\end{equation*}

Dividing by $\|Z\|^{p-1}_{L^p}$ and integrating in time yield
\begin{equation}
\label{eq:3.4}\begin{split}
\|Z\|_{L^p}\leq \|Z_0\|_{L^p}+C \int_{0}^{t}\big(\|\omega(.,\tau)\|_{L^p}+\|\nabla\theta(.,\tau)\|_{L^p}\big)d\tau.
\end{split}
\end{equation}

In particular,
\begin{equation}
\label{eq3.5}\begin{split}
\|Z\|_{L^\infty}\leq \|Z_0\|_{L^\infty}+C \int_{0}^{t}\big(\|\omega(.,\tau)\|_{L^\infty}+\|\nabla\theta(.,\tau)\|_{L^\infty}\big)d\tau.
\end{split}
\end{equation}
From the second equation of \eqref{equ:1.1}, we have

\begin{equation*}
\label{3.5}\begin{split}
\|\omega\|_{L^p}\leq \|\omega_0\|_{L^p}+C \int_{0}^{t}\|\nabla\times u(.,\tau)\|_{L^p}d\tau,
\end{split}
\end{equation*}
which together with \eqref{eq3.2} implies that
\begin{equation}
\label{3.6}\begin{split}
\|\omega\|_{L^p}\leq \|\omega_0\|_{L^p}+C \int_{0}^{t}\big(\|\omega(.,\tau)\|_{L^p}+\|Z(.,\tau)\|_{L^p}\big)d\tau,
\end{split}
\end{equation}

and
taking limit $p\to\infty $ lead to
\begin{equation}
\label{3.7}\begin{split}
\|\omega\|_{L^\infty}&\leq \|\omega_0\|_{L^\infty}+C \int_{0}^{t}\big(\|\omega(.,\tau)\|_{L^\infty}+\|Z(.,\tau)\|_{L^\infty}\big)d\tau.
\end{split}
\end{equation}

On the other hand, we need make  $L^{p}$-estimate of $\nabla\theta$. Taking $\nabla^{\bot} = (-\partial_{x2}, \partial_{x1} )$ to the third equation of \eqref{equ:1.1}, we obtain
\begin{equation}
\label{eq3.8} \nabla^{\bot}\theta_t + u\cdot \nabla \nabla^{\bot}\theta-\mu\Delta \nabla^{\bot}\theta =\nabla^{\bot}\theta\cdot\nabla u+ \nabla u.
\end{equation}
Taking $L^{2}$ inner product \eqref{eq3.8} with  $| \nabla^{\bot}\theta|^{p-2} \nabla^{\bot}\theta$$(2\leq p<\infty)$,  using the incompressibility of $u$ and the  H\"{o}lder inequality, we have
\begin{equation*}
\label{3.25}\begin{split}
&\frac{1}{p}\frac{d}{dt}\|\nabla\theta\|^p_{L^p}\leq \|\nabla u\|_{L^p}\|\nabla\theta\|_{L^\infty}\|\nabla\theta\|^{p-1}_{L^p}+\|\nabla u\|_{L^p}\|\nabla\theta\|^{p-1}_{L^p}.
\end{split}
\end{equation*}

Dividing by $\|\nabla\theta\|^{p-1}_{L^p}$ and integrating in time yield
\begin{equation}
\label{eq3.9}\begin{split}
\|\nabla\theta\|_{L^p}&\leq \|\nabla\theta_0\|_{L^p}+C \int_{0}^{t}\big(1+\|\nabla\theta(.,\tau)\|_{L^\infty}\big)\|\nabla u\|_{L^p}d\tau\\
&\leq \|\nabla\theta_0\|_{L^p}+C \int_{0}^{t}\big(1+\|\nabla\theta(.,\tau)\|_{L^\infty}\big)\|\nabla\times u\|_{L^p}d\tau\\
&\leq \|\nabla\theta_0\|_{L^p}+C \int_{0}^{t}\big(1+\|\nabla\theta(.,\tau)\|_{L^\infty}\big)\big(\|Z\|_{L^p}+\|\omega\|_{L^p}\big)d\tau.
\end{split}
\end{equation}

Due to \begin{equation}
\label{eq3.91}\|\nabla \theta(.,\tau)\|_{L^{\infty}}< +\infty,\end{equation}{combining with \eqref{eq:3.4}, \eqref{3.6} and \eqref{eq3.9}, and applying the Gronwall inequality, we have
	\begin{equation}
	\label{eq3.10}\begin{split}
	\|Z\|_{L^p}+\|\omega\|_{L^p}+\|\nabla\theta\|_{L^p}&\leq \Big(\|\Omega_0\|_{L^p}+\|\omega_0\|_{L^p}+\|\nabla\theta_0\|_{L^p}\Big)\exp\Big( C \int_{0}^{t}\big(1+\|\nabla\theta(.,\tau)\|_{L^\infty}\big)d\tau\Big),
	\end{split}
	\end{equation}
	for $2\leq p<\infty.$
	
Furthermore, employing \eqref{eq3.5} and \eqref{3.7}, and then applying the Gronwall inequality, we obtain
\begin{equation}
\label{eq3.11}\begin{split}
\|Z\|_{L^\infty}+\|\omega\|_{L^\infty}&\leq C\big(\|\Omega_0\|_{L^\infty}+\|\omega_0\|_{L^\infty}+\|\nabla\theta_0\|_{L^p}\big)e^{Ct }\leq C(u_0,\omega_0,\theta_0;T).
\end{split}
\end{equation}
From (4.12) and  \eqref{eq3.11}, we finally conclude that
\begin{equation}
\label{eq3.12}\begin{split}
\|\Omega\|_{L^p}+\|\Omega\|_{L^\infty}<+\infty.
\end{split}
\end{equation}
Now, recall the following well-known result (see e.g. \cite{BKM}):

\begin{equation*}
\label{eq3.13}\begin{split}
\|\nabla u\|_{L^\infty}\leq C\Big\{1+\|\Omega\|_{L^p}+\|\Omega\|_{L^\infty}\ln(1+\|u\|_{H^{s}})\Big\}.
\end{split}
\end{equation*}
From \eqref{eq3.12}, the above inequality  implies that
\begin{equation}
\label{eq3.14}\begin{split}
\|\nabla u\|_{L^\infty}\leq C\ln\big(1+\|u\|_{H^{s}}\big).
\end{split}
\end{equation}

We write the  equation of $\omega(t)$ in  \eqref{equ:1.1} as
\begin{equation}\label{eq4.17}
\begin{array}{ll}
\omega(t)=e^{\gamma t \Delta}u_0+\displaystyle\int_{0}^{t}e^{\gamma(t-\tau)\Delta}(2\kappa \nabla\times\, u-4\kappa\omega-u\cdot\nabla\omega)d\tau.
\end{array}
\end{equation}
Taking  $p=2$ and $2\leq q < \infty$  in Lemma 2.5, we get
\begin{equation}\label{eq4.18}
\begin{array}{ll}
\displaystyle \int^t_{0} \|\Delta \omega\|^2_{L^q}d\tau \leq C\|u_0\|^2_{H^2} +C\displaystyle \int^t_{0} \|\nabla\times\, u\|^2_{L^q}+\|\omega\|^2_{L^q}+\|u\cdot\nabla\omega\|^2_{L^q}d\tau.
\end{array}
\end{equation}

On the other hand, using the Sobolev embedding and combining with (4.1) and (4.2), we have
\begin{equation}\label{eq4.19}
\begin{array}{ll}
\displaystyle \int^t_{0}\|\omega\|^2_{L^q}d\tau \leq C \displaystyle \int^t_{0}\|\nabla \omega\|^2_{L^2}d\tau\leq C(u_0,\omega_0,\theta_0;T),
\end{array}	
\end{equation}
\begin{equation}\label{eq4.20}
\begin{array}{ll}
\displaystyle \int^t_{0}\|u\cdot\nabla \omega\|^2_{L^q}d\tau&\leq\displaystyle \int^t_{0}\|u\|^2_{L^{2q}}\|\nabla \omega\|^2_{L^{2q}}d\tau
\\&\leq C\displaystyle \int^t_{0}\|\nabla u\|^2_{L^{2}}\|\nabla^2 \omega\|^2_{L^{2}}d\tau\\
&\leq C \sup\limits_{t\in(0,T)}\|\nabla u\|^2_{L^2}\displaystyle\int_{0}^{t}\|\nabla^2\omega\|^2_{L^2}d\tau\\&\leq
C(u_0,\theta_0,\omega_0;T).
\end{array}	
\end{equation}

Thus, from \eqref{eq4.18}  we get
\begin{equation}\label{eq4.21}
\begin{array}{ll}
\displaystyle \int^t_{0}\|\Delta \omega\|^2_{L^q}d\tau\leq C\displaystyle \int^t_{0} \|\nabla\times\, u\|^2_{L^q} +C(u_0,\omega_0,\theta_0;T).
\end{array}	
\end{equation}
Employing (4.14) and (4.21) and then using  the Sobolev embedding, we obtain
\begin{equation}\label{eq4.22}
\displaystyle \int_{0}^{t}\|\nabla \omega\|_{L^{\infty}}\leq C(u_0,\omega_0,\theta_0;T).
\end{equation}
Exactly as in the proof of \eqref{eq 3.5}, we also have

\begin{equation}
\label{eq 3.15}\begin{split}\ln\big(\|u\|_{H^{s}}^{2}
+\|\omega\|_{H^{s}}^{2}+\|\theta\|_{H^{s}}^{2}\big)
&\leq C\ln \big(\|u_0\|_{H^{s}}^{2}
+\|\omega_0\|_{H^{s}}^{2}+\|\theta_0\|_{H^{s}}^{2}\big)\\&\quad+C
\int_{0}^{t}\big(1+ \|\nabla u\|_{L^{\infty}}+\|\nabla\omega\|_{L^{\infty}}+\|\nabla
\theta\|_{L^{\infty}}\big)d\tau.
\end{split}
\end{equation}

 Combining with \eqref{eq3.91}, \eqref{eq3.14}, (4.22) and \eqref{eq 3.15}, and then applying the Gronwall inequality,  we conclude that
\begin{equation}
\|u\|^2_{H^s}+\|\omega\|^2_{H^s}+\|\theta\|^2_{H^s}\leq C(u_0,\omega_0,\theta_0;T).
\end{equation}
The proof of Theorem \ref{th2} is completed.
\section{The  proof of Theorem \ref{th3}}
In order to prove  the global regularity
part of Theorems \ref{th3}, it suffices to get estimate \eqref{blow-up} for all $T\in (0,\infty)$
for the smooth solutions of the  system \eqref{equ:1.1}. Here, we
divide it into  the following two parts.
\subsection{ Bound for $\|u\|_{L^{\infty}}$}
\ \ \ \ \

Now, we recall the equation about the $\Omega=\nabla\times u$
\begin{equation}\label{eq5.3}
\Omega_t + u\cdot\nabla \Omega +2\kappa \Delta  \omega =\partial_1\theta.
\end{equation}

Taking scaler product \eqref{eq5.3} in $L^2$ by $|\Omega|^{p-2}\Omega(2<p<\infty)$, after integration by part, and using the incompressibility of $u$, we obtain the inequality
$$\frac{1}{p}\frac{d}{dt}\|\Omega\|^{p}_{L^p}\leq \|\Delta \omega\|_{L^p}\|\Omega\|^{p-1}_{L^p}+\|\nabla\theta\|_{L^p}\|\Omega\|^{p-1}_{L^p}.$$
Dividing by $\|\nabla\theta\|^{p-1}_{L^p}$, integrating in time  and then  using the Sobolev embedding yield
\begin{equation}\label{eq5.4}
\begin{array}{ll}
\|\Omega\|_{L^p}&\leq \|\Omega_0\|_{L^p} +\displaystyle \int^t_{0}\|\Delta \omega\|_{L^p}d\tau +\displaystyle \int^t_{0} \|\nabla\theta\|_{L^p} d\tau\\
&\leq \|\Omega_0\|_{L^p}+\displaystyle \int^t_{0}\|\Delta \omega\|_{L^p}d\tau +\displaystyle \int^t_{0} \|\nabla^2\theta\|_{L^2} d\tau \\
&\leq\|\Omega_0\|_{L^p}+C(T)  \Big(\displaystyle \int^t_{0}\|\Delta \omega\|^2_{L^p}d\tau\Big)^{\frac{1}{2}} +\Big(\displaystyle \int^t_{0} \|\nabla^2\theta\|^2_{L^2} d\tau\Big)^{\frac{1}{2}},
\end{array}
\end{equation}

which implies  that $$\|\Omega\|^2_{L^p}\leq\|\Omega_0\|^2_{L^p}+C(T)\Big(\displaystyle \int^t_{0}\|\Delta \omega\|^2_{L^p}d\tau +\displaystyle \int^t_{0} \|\nabla^2\theta\|^2_{L^2} d\tau\Big). $$
Using the Gronwall inequality and combining with (4.2) and (4.21), we finally get
\begin{equation}\label{eq5.5}
\|\Omega\|_{L^p}\leq C(u_0,\omega_0,\theta_0;T).
\end{equation}

Therefore, by the following Gagliardo-Nirenberg interpolation inequality in $\mathbb{R}^{2}$
\begin{equation}
\label{Gagliardo-Nirenberg interpolation}\|f\|_{L^{\infty}}\leq C\|f\|^{\frac{p-2}{2p-2}}_{L^{2}}\|\nabla f\|^{\frac{p}{2p-2}}_{L^{p}},\quad f\in W^{1,p}, \quad p>2,\end{equation}
and the Calderon-Zygmund inequality and  \eqref{eq4.1}, \eqref{eq5.5},  we have
\begin{equation}\label{eq5.7}
\|u\|_{L^{\infty}}\leq C\|u\|^{\frac{p-2}{2p-2}}_{L^{2}}\|\nabla u\|^{\frac{p}{2p-2}}_{L^{p}}\leq C\|u\|^{\frac{p-2}{2p-2}}_{L^{2}}\|\Omega\|^{\frac{p}{2p-2}}_{L^{p}}\leq C\big(u_0,\omega_0,\theta_0,T,p\big).\end{equation}
\subsection{ Bound for $\|\nabla \theta\|_{L^{\infty}}$}
\ \ \ \ \ 

Taking scalar product \eqref{eq3.8} in $L^{2}$ by $| \nabla^{\bot}\theta|^{p-2} \nabla^{\bot}\theta$ $(2< p<\infty)$ again, after integration by
part, and using the incompressibility of $u$,  the   H\"{o}lder inequality, the  Young inequality and the Calderon-Zygmund inequality, we get

\begin{equation*}
\label{3.25}\begin{split}
&\frac{1}{p}\frac{d}{dt}\|\nabla^{\bot}\theta\|^p_{L^p}+\mu(p-1)\int|D^{2} \theta|^{2}| \nabla^{\bot}\theta|^{p-2}dx\\
&=\int\nabla^{\bot}\theta\cdot\nabla u|\nabla^{\bot}\theta|^{p-2} \nabla^{\bot}\theta dx+\int\nabla u|\nabla^{\bot}\theta|^{p-2} \nabla^{\bot}\theta dx\\
&\leq \int |u|^{2}|\nabla^{\bot}\theta|^{p}dx+\frac{\mu(p-1)}{2}\int |\nabla^{2}\theta|^{2}\nabla^{\bot}\theta|^{p-2}dx+\|\nabla u\|_{L^p}\|\nabla\theta\|^{p-1}_{L^p}\\
&\leq C\|u\|^{2}_{L^\infty}\|\nabla\theta\|^{p}_{L^p}+\frac{\mu(p-1)}{2}\int |\nabla^{2}\theta|^{2}\nabla^{\bot}\theta|^{p-2}dx+\|\Omega\|^{p}_{L^p}+ \|\nabla\theta\|^{p}_{L^p},
\end{split}
\end{equation*}
which implies  that
\begin{equation*}
\label{3.25}\begin{split}
\frac{1}{p}\frac{d}{dt}\|\nabla^{\bot}\theta\|^p_{L^p}+\frac{\mu(p-1)}{2}\int|D^{2} \theta|^{2}| \nabla^{\bot}\theta|^{p-2}dx\leq C\big(1+\|u\|^{2}_{L^\infty}\big)\|\nabla\theta\|^{p}_{L^p}+\|\Omega\|^{p}_{L^p},
\end{split}
\end{equation*}

which together with \eqref{eq5.5}, \eqref{eq5.7} and  the Gronwall inequality implies that
\begin{equation} \label{eq5.9}
\begin{split}
\|\nabla \theta\|_{L^p}\leq C\big(u_0,\omega_0,\theta_0,T,p\big).
\end{split}
\end{equation}
Taking operation $D^{2}$ on  the third equation of \eqref{equ:1.1}, and  then taking  $L^{2}$ scalar product of this  with $| D^{2}\theta|^{p-2} D^{2}\theta$ $(2< p<\infty)$, after integration by
part, and  using the    H\"{o}lder inequality, the Young inequality, Gagliardo-Nirenberg interpolation inequality \eqref{Gagliardo-Nirenberg interpolation}  and the Calderon-Zygmund inequality, we deduce that
\begin{equation*}
\label{3.25}\begin{split}
\frac{1}{p}\frac{d}{dt}&\|D^{2}\theta\|^p_{L^p}+\mu(p-1)\int|D^{3} \theta|^{2}| D^{2}\theta|^{p-2}dx\\
&=-\int D^{2}(u\cdot\nabla \theta)|D^{2}\theta|^{p-2} D^{2}\theta dx+\int D^{2}u|D^{2}\theta|^{p-2} D^{2}\theta dx\\
&=(p-1)\int Du\cdot\nabla \theta|D^{2}\theta|^{p-2} D^{3}\theta dx+(p-1)\int u\cdot D\nabla\theta|D^{2}\theta|^{p-2} D^{3}\theta dx
\\&\quad+(p-1)\int Du|D^{2}\theta|^{p-2} D^{3}\theta dx\\
&\leq  C\|\nabla\theta\|^{2}_{L^{\infty}}\|D u\|^{2}_{L^p}\|D^{2}\theta\|^{p-2}_{L^p}+\frac{\mu(p-1)}{2}\int|D^{3} \theta|^{2}| D^{2}\theta|^{p-2}dx+C\|u\|^{2}_{L^{\infty}}\|D^{2}\theta\|^{p}_{L^p}\\&\quad+\|D u\|^{2}_{L^p}\|D^{2}\theta\|^{p-2}_{L^p}\\
&\leq  C\|\nabla\theta\|^{2-\frac{4}{p}}_{L^{p}}\|\Omega\|^{2}_{L^p}\|D^{2}\theta\|^{\frac{4}{p}+p-2}_{L^p}+\frac{\mu(p-1)}{2}\int|D^{3} \theta|^{2}| D^{2}\theta|^{p-2}dx+C\|u\|^{2}_{L^{\infty}}\|D^{2}\theta\|^{p}_{L^p}\\&\quad+\|\Omega\|^{2}_{L^p}\|D^{2}\theta\|^{p-2}_{L^p},
\end{split}
\end{equation*}

where we have used the following  interpolation inequality in $\mathbb{R}^{2}$
\begin{equation}\label{5.10}
\|f\|_{L^{\infty}}\leq C\|f\|^{1-\frac{2}{p}}_{L^{p}}\|\nabla f\|^{\frac{2}{p}}_{L^{p}},\quad f\in W^{1,p}, \quad p>2.\end{equation}
Note $\frac{4}{p}+p-2< p$ when $p > 2$, thanks to the estimates \eqref{eq5.5} \eqref{eq5.7} and \eqref{eq5.9}, we obtain
\begin{equation*}
\label{3.25}\begin{split}
\frac{d}{dt}\|D^{2}\theta\|^p_{L^p}\leq C+C\|D^{2}\theta\|^p_{L^p}.
\end{split}
\end{equation*}
Using the Gronwall
inequality yields that
\begin{equation}\label{eq5.11}
\begin{split}
\|D^{2}\theta\|_{L^p}\leq C\big(u_0,\omega_0,\theta_0,T,p\big).
\end{split}
\end{equation}
Furthermore, employing (5.9) and \eqref{eq5.11}, and using Gagliardo-Nirenberg interpolation \eqref{5.10}, we finally deduce that
\begin{equation*}
\begin{split}
\|\nabla\theta\|_{L^\infty}\leq C(u_0,\omega_0,\theta_0;T) \quad \forall t\in [0,T],
\end{split}
\end{equation*}
where $C=C(\|u_0\|_{W^{2,p}},\|\omega_0\|_{W^{2,p}},\|\theta_0\|_{W^{2,p}})$. Due to the embedding $H^{s}(\mathbb{R}^{2})\hookrightarrow W^{2,p}(\mathbb{R}^{2})$, for all $s> 2$ and $p > 2$, and thus we attained
estimate \eqref{blow-up} for all $T\in (0,\infty)$ and for all $u_0, \omega_0, \theta_0 \in H^{s}(\mathbb{R}^{2})$ with $s>2$.
The proof of Theorem \ref{th3} is completed.

{\bf acknowledgements:} {\rm We are grateful to Professor  Yi  Zhou for his guidance and encouragement, which greatly improved our original manuscript. }

\begin{center}

\end{center}

\end{document}